\documentstyle[%
  12pt,%
  amscd%
  ]{amsart}

\theoremstyle{plain}
  \newtheorem{thm}{Theorem}[section]
  \newtheorem{prop}[thm]{Proposition}
  \newtheorem{lem}[thm]{Lemma}
  \newtheorem{cor}[thm]{Corollary}
\theoremstyle{definition}
  \newtheorem{dfn}[thm]{Definition}
\numberwithin{equation}{thm}

\newcommand{\spec}{\operatorname{Spec}}
\newcommand{\proj}{\operatorname{Proj}}
\newcommand{\ext}{\operatorname{Ext}}
\renewcommand{\hom}{\operatorname{Hom}}
\newcommand{\ann}{\operatorname{ann}}
\newcommand{\depth}{\operatorname{depth}}
\newcommand{\coker}{\operatorname{Coker}}
\renewcommand{\ker}{\operatorname{Ker}}
\newcommand{\qtn}{\mathop{:}}                
\newcommand{\angled}[1]{\langle #1 \rangle}  
\newcommand{\trans}[1]{\overline{#1}}
\newcommand{\closure}{^*}  
\newcommand{\hlz}[1]{_{(#1)}}

\newcommand{\dlim}{\mathop{\varinjlim}}
\newcommand{\dc}{D^\bullet}  

\newcommand{\sheafO}{{\cal O}}
\newcommand{\ideala}{{\frak a}}
\newcommand{\idealb}{{\frak b}}
\newcommand{\idealc}{{\frak c}}
\newcommand{\ideall}{{\frak l}}
\newcommand{\idealm}{{\frak m}}
\newcommand{\idealM}{{\frak M}}
\newcommand{\idealn}{{\frak n}}
\newcommand{\idealN}{{\frak N}}
\newcommand{\idealp}{{\frak p}}
\newcommand{\idealq}{{\frak q}}
\newcommand{\idealr}{{\frak r}}

\begin{document}


  \title[Macaulayfication]{%
    On Macaulayfication of certain \\
      quasi-projective schemes}

  \author[T.~Kawasaki]{%
    Takesi Kawasaki}

  \address{%
    Department of Mathematics,
    Tokyo Metropolitan University,
    Hachioji Minami-Ohsawa 1--1,
    Tokyo 192--03 Japan}

  \email{kawasaki@@math.metro-u.ac.jp}

  \thanks{%
    The author is partially supported
    by Grant-Aid for Co-operative Research.}



\maketitle

\begin{abstract}

  The Macaulayfication of a Noetherian scheme~$X$ is
    a birational proper morphism~$Y \rightarrow X$
    such that $Y$ is a Cohen-Macaulay scheme.
  Of course,
    a desingularization is a Macaulayfication
    and Hironaka gave a desingularization
    of arbitrary algebraic variety
    over a field of characteristic~$0$.
  In~1978
    Faltings gave a Macaulayfication of a quasi-projective scheme
    whose non-Cohen-Macaulay locus is of dimension $0$~or~$1$
    by a characteristic-free method.

  In the present article,
    we shall construct a Macaulayfication
    of a quasi-projective scheme
    if the dimension of its non-Cohen-Macaulay locus
    is at most~$2$.
  Of course,
    our method is independent of the characteristic of a scheme.

\end{abstract}

\section{Introduction}

  Let $X$ be a Noetherian scheme.
  A birational proper morphism~$Y \rightarrow X$ of schemes
    is said to be
    a {\em Macaulayfication\/} of~$X$
    if $Y$ is a Cohen-Macaulay scheme.
  This notion was introduced by Faltings~\cite{Faltings:78:Macaulay}
    and he established that
    there exists a Macaulayfication of a quasi-projective scheme
    over a Noetherian ring possessing a dualizing complex
    if its non-Cohen-Macaulay locus
    is of dimension $0$~or~$1$.
  Of course,
    a desingularization is a Macaulayfication
    and Hironaka gave a desingularization
    of arbitrary algebraic variety over a field of characteristic~$0$.
  But Faltings' method to construct a Macaulayfication
    is independent of the characteristic of a scheme.
  Furthermore,
    several authors are interested in a Macaulayfication.

  For example,
    Goto and Schenzel independently showed
    the converse of Faltings' result in a sense.
  Let $A$ be a Noetherian local ring possessing a dualizing complex,
    hence its non-Cohen-Macaulay locus is closed,
    and assume that $\dim A / \idealp = \dim A$
    for any associated prime ideal~$\idealp$ of~$A$.
  Then the non-Cohen-Macaulay locus of~$A$ consists of
    only the maximal ideal
    if and only if
    $A$ is a generalized Cohen-Macaulay ring
    but not a Cohen-Macaulay ring~\cite{Schenzel:75:einige}.
  When this is the case,
    Faltings~\cite[Satz 2]{Faltings:78:Macaulay} showed that
    there exists a parameter ideal~$\idealq$ of~$A$
    such that the blowing-up~$\proj A[\idealq t]$ of~$\spec A$
    with center~$\idealq$ is Cohen-Macaulay,
    where $t$ denotes an indeterminate.
  Conversely,
    Goto~\cite{Goto:82:blowing} proved that
    if there is a parameter ideal~$\idealq$ of~$A$
    such that $\proj A[\idealq t]$ is Cohen-Macaulay,
    then $A$ is a generalized Cohen-Macaulay ring.
  Moreover,
    he showed that $A$ is Buchsbaum if and only if
    $\proj A[\idealq t]$ is Cohen-Macaulay
    for every parameter ideal~$\idealq$ of~$A$:
    see also \cite{Schenzel:83:standard}.

  Brodmann~\cite{Brodmann:83:local} also studied
    the blowing-up of a generalized Cohen-Macaulay ring
    with center a parameter ideal.
  Furthermore,
    he constructed Macaulayfications
    in a quite different way from Faltings.
  Let $A$ be a Noetherian local ring possessing a dualizing complex.
  We let $d= \dim A$
    and $s$ be the dimension of its non-Cohen-Macaulay locus.
  If $s=0$,
    then Brodmann~\cite[Proposition 2.13]{Brodmann:83:two}
    gave an ideal~$\idealb$ of height~$d-1$
    such that $\proj A[\idealb t]$ is Cohen-Macaulay.
  If $s=1$, then
    Faltings' Macaulayfication~\cite[Satz 3]{Faltings:78:Macaulay}
    of~$\spec A$ consists of
    two consecutive blowing-ups $Y \rightarrow X \rightarrow \spec A$
    where the center of the first blowing-up
    is an ideal of height~$d-1$.
  In this case,
    Brodmann gave two other Macaulayfications of~$\spec A$:
    the first one~\cite{Brodmann:80:Kohomologische}
      is the composite of a blowing-up $X \rightarrow \spec A$
      with center an ideal of height~$d-1$
      and a finite morphism~$Y \rightarrow X$;
    the second one~\cite[Corollary 3.11]{Brodmann:83:two}
      consists of
      two consecutive blowing-ups $Y \rightarrow X \rightarrow \spec A$
      where the center of the first blowing-up is
      an ideal of height~$d-2$.

  In this article,
    we are interested in a Macaulayfication of
    the Noetherian scheme
    whose non-Cohen-Macaulay locus is of dimension~$2$.
  Let $A$ be a Noetherian ring possessing a dualizing complex
    and $X$ a quasi-projective scheme over~$A$.
  Then $X$ has a dualizing complex with codimension function~$v$.
  Furthermore the non-Cohen-Macaulay locus~$V$ of~$X$
    is closed.
  We define a function $u \colon X \to {\Bbb Z}$
    to be $u(p) = v(p) + \dim \overline{\{p\}}$.
  We will establish the following theorem:

  \begin{thm} \label{mthm}
    If $\dim V \leq 2$ and $u$ is locally constant on~$V$,
      then $X$ has a Macaulayfication.
  \end{thm}

  If $\dim V \leq 1$,
    then $u$ is always locally constant on~$V$.
  Therefore, this theorem contains Faltings' result.
  Furthermore,
    we note if $X$ is a projective scheme
    over a Gorenstein local ring,
    then $u$ is constant on~$X$.

  We agree that
    $A$ denotes a Noetherian local ring with maximal ideal~$\idealm$
    except for Section~\ref{sec:6}.
    Assume that $d = \dim A > 0$.
  We refer the reader to~%
      \cite{Hartshorne:66:residue,%
      Hartshorne:77:algebraic,%
      Matsumura:89:commutative,%
      Stuckrad-Vogel:86:Buchsbaum}
    for unexplained terminology.

\section{Preliminaries}

  In this section,
    we state some definitions and properties
    of a local cohomology and an ideal transform.
  Let $\idealb$ be an ideal of~$A$.

  \begin{dfn}
    The local cohomology functor~$H_\idealb^p(-)$ and
      the ideal transform functor~$D_\idealb^p(-)$
      with respect to~$\idealb$
      are defined to be
        $$
          H_\idealb^p(-) = \dlim_m \ext_A^p(A/\idealb^m, -)
            \quad
            \text{and}
            \quad
          D_\idealb^p(-) = \dlim_m \ext_A^p(\idealb^m, -),
        $$
      respectively.
  \end{dfn}

  For an $A$-module $M$,
    there exists an exact sequence
      \begin{equation} \label{eqn:2.1.1}
        0 @>>>
        H_\idealb^0(M) @>>>
        M @>\iota>>
        D_\idealb^0(M) @>>>
        H_\idealb^1(M) @>>>
        0
      \end{equation}
    and isomorphisms
      $$
        D_\idealb^p(M) \cong H_\idealb^{p+1}(M)
        \quad
        \text{for all $p>0$}.
      $$
  They induces that
      \begin{equation} \label{eqn:2.1.2}
        H_\idealb^p D_\idealb^0(M) =
        \begin{cases}
          0, & p=0, 1;
        \\
          H_\idealb^p(M), & \text{otherwise}.
        \end{cases}
      \end{equation}
  If $\idealb$ contains an $M$-regular element~$a$,
    then we can regard $D_\idealb^0(M)$ as a submodule
    of the localization~$M_a$ with respect to~$a$
    and $\iota$ is the inclusion.

  It is well-known that $H_\idealb^p(-)$ is naturally isomorphic to
    the direct limit of Koszul cohomology.
  In particular,
    let $\idealb = (f_1, \dots, f_h)$ and
    $M$ be an $A$-module.
  Then
      $$
        H_\idealb^h(M) = \dlim_m M/(f_1^m, \dots, f_h^m)M
          \quad
          \text{and}
          \quad
        H_\idealb^0(M) = \bigcap_{i=1}^h 0 \qtn_M \angled{f_i},
      $$
    where $0 \qtn \angled{f_i}$ denotes
    $\bigcup_{m = 1}^\infty 0 \qtn f_i^m$.
  Furthermore,
    let $A \rightarrow B$ be a ring homomorphism.
  Then there exists a natural isomorphism
    $H_\idealb^p(M) \cong H_{\idealb B}^p(M)$
    for a $B$-module~$M$.

  The following lemma is frequently used in this article.

  \begin{lem}[Brodmann \cite{Brodmann:83:Einige}]
      \label{lem:2.2}
    Let $\idealb = (f_1, \dots, f_h)$ and
      $\idealc = (f_1, \dots, f_{h-1})$
      be two ideals.
    Then there exists a natural long exact sequence
        $$
          \cdots @>>>
          [H_\idealc^{p-1}(-)]_{f_h} @>>>
          H_\idealb^p(-) @>>>
          H_\idealc^p(-) @>>>
          [H_\idealc^p(-)]_{f_h} @>>>
          \cdots.
        $$
  \end{lem}

  Next we state on the annihilator of local cohomology modules.

  \begin{dfn}
    For any finitely generated $A$-module~$M$,
      we define an ideal~$\ideala_A(M)$ to be
        $$
          \ideala_A(M) = \prod_{p=0}^{\dim M -1}
          \ann H_\idealm^p(M).
        $$
  \end{dfn}

  We note that a finitely generated $A$-module~$M$ is Cohen-Macaulay
    if and only if $\ideala_A(M) = A$,
    and that $M$ is generalized Cohen-Macaulay
    if and only if $\ideala_A(M)$ is an $\idealm$-primary ideal.
  The notion of~$\ideala_A(-)$ plays a key role in this article.
  In fact,
    Schenzel~\cite{Schenzel:79:dualizing} showed that
    $V(\ideala_A(A))$ coincides with the non-Cohen-Macaulay locus of~$A$
    if it possesses a dualizing complex and is equidimensional.
  He also gave the following lemma~%
    \cite{Schenzel:79:dualizing,Schenzel:82:cohomological}:

  \begin{lem} \label{lem:2.4}
    Let $M$ be a finitely generated $A$-module and
      $x_1$,~\dots, $x_n$ a system of parameters for~$M$.
    Then
        $
          (x_1, \dots, x_{i-1})M \qtn x_i \subseteq
          (x_1, \dots, x_{i-1})M \qtn \ideala_A(M)
        $
      for any $1 \leq i \leq n$.
    In particular,
      if $x_i \in \ideala_A(M)$,
      then the equality holds.
  \end{lem}

  Let $R = \bigoplus_{n \geq 0} R_n$ be a Noetherian graded ring
    where $R_0 = A$.
  A graded module~$M = \bigoplus_n M_n$ is said to be
    {\em finitely graded\/}
    if $M_n=0$ for all but finitely many~$n$.
  The following lemma is an easy consequence of~\cite{Faltings:78:uber}.

  \begin{lem} \label{lem:2.5}
    Let $\idealb$ be a homogeneous ideal of~$R$
      containing~$R_+ = \bigoplus_{n>0} R_n$
      and $M$ a finitely generated graded $R$-module.
    We assume that $A$ possesses a dualizing complex.
    Let $p$ be the largest integer
      such that, for all~$q \leq p$,
      $H_\idealb^q(M)$ is finitely graded.
    Then $\depth M\hlz\idealp \geq p$
      for any closed point~$\idealp$ of~$\proj R$,
      that is, $\idealp$ is a homogeneous prime ideal
      such that $\dim R/ \idealp = 1$
      and $R_+ \not\subseteq \idealp$.
  \end{lem}

\section{A Rees algebra obtained by an ideal transform}
  \label{sec:3}

  \begin{dfn}
    A sequence $f_1$,~\dots, $f_h$ of elements of~$A$
      is said to be a d-sequence on an $A$-module~$M$ if
        $
          (f_1, \dots, f_{i-1})M \qtn f_i f_j
          = (f_1, \dots, f_{i-1})M \qtn f_j
        $
      for any $1 \leq i \leq j \leq h$.

    We shall say that $f_1$,~\dots, $f_h$ is
      an unconditioned strong d-sequence
      (for short, {\em u.s.d-sequence\/})
      on~$M$
      if $f_1^{n_1}$,~\dots, $f_h^{n_h}$ is a d-sequence on~$M$
      in any order and
      for arbitrary positive integers $n_1$,~\dots, $n_h$.
  \end{dfn}

  The notion of u.s.d-sequences was introduced
    by Goto and Yamagishi~\cite{Goto-Yamagishi::theory}
    to refine arguments on Buchsbaum rings
    and generalized Cohen-Macaulay rings.
  Their theory contains
    Brodmann's study on the Rees algebra with respect to an ideal
    generated by a pS-sequences~\cite{Brodmann:83:local}.
  But Brodmann~\cite{Brodmann:84:local} also studied
    the ideal transform of such a Rees algebra.
  The purpose of this section is to study
    an ideal transform of the Rees algebra with respect to an ideal
    generated by a u.s.d-sequence.

  Let $f_0$,~\dots, $f_h$ be a sequence of elements of~$A$
    where $h \geq 1$ and
    $\idealq = (f_1, \dots, f_h)$.

  \begin{lem} \label{lem:3.2}
    If $f_1$,~\dots, $f_h$ be a d-sequence on~$A/ f_0A$,
      then
        $$
          [(f_1, \dots, f_k) \idealq^n] \qtn f_0 =
          (f_1, \dots, f_k) [\idealq^n \qtn f_0] +
          0 \qtn f_0
        $$
      for any $1 \leq k \leq h$ and $n>0$.
  \end{lem}

  \begin{pf}
    It is obvious that
      the left hand side contains the right one.
    We shall prove the inverse inclusion
      by induction on~$k$.
    Let $a$ be an element of the left hand side.

    When $k=1$,
      we put $f_0 a = f_1 b$ where $b \in \idealq^n$.
    By using \cite[Theorem 1.3]{Goto-Yamagishi::theory},
      we obtain
      $b \in (f_0) \qtn f_1 \cap \idealq^n \subseteq (f_0)$.
    If we put $b = f_0 a'$,
      then $a' \in \idealq^n \qtn f_0$ and
      $f_0 (a - f_1 a') =0$.
    Thus we get
      $a \in f_1 [\idealq^n \qtn f_0] + 0 \qtn f_0$.

    When $k > 1$,
      we put $f_0 a = b + f_k c$
      where
        $
          b \in (f_1, \dots, f_{k-1}) \idealq^n
        $
      and $c \in \idealq^n$.
    Then we obtain
        \begin{align*}
          c & \in (f_0, \dots, f_{k-1}) \qtn f_k \cap \idealq^n
        \\
          & \subseteq (f_0) + (f_1, \dots, f_{k-1}) \idealq^{n-1}
        \end{align*}
      by using~\cite[Theorem 1.3]{Goto-Yamagishi::theory} again.
    If we put $c = f_0 a' + b'$
      where
        $$
          b' \in (f_1, \dots, f_{k-1}) \idealq^{n-1},
        $$
      then $a' \in \idealq^n \qtn f_0$.
    Thus we get
        \begin{align*}
          a - f_k a' & \in [(f_1, \dots, f_{k-1}) \idealq^n] \qtn f_0
        \\
          & = (f_1, \dots, f_{k-1}) [\idealq^n \qtn f_0] + 0 \qtn f_0
        \end{align*}
      by induction hypothesis.
    The proof is completed.
  \end{pf}

  Let $\trans\idealq = \idealq \qtn \angled{f_0}$.
  If $f_0$ is $A$-regular
    and $f_1$,~\dots, $f_h$ is a d-sequence
    on~$A/ f_0^lA$ for all~$l > 0$,
    then Lemma~\ref{lem:3.2} assures us that
      \begin{equation} \label{eqn:3.2.1}
        \idealq^{n-1} \trans\idealq =
        \trans\idealq^n =
        \idealq^n \qtn \angled{f_0}
          \quad
          \text{for all $n>0$}.
      \end{equation}
  Therefore the Rees algebra $\trans R = A[\trans\idealq t]$
    is finitely generated over $R = A [\idealq t]$.
  The following is an analogue of~\cite[Lemma 3.4]{Goto:82:blowing}.

  \begin{thm} \label{thm:3.3}
    Let $B = A[\trans\idealq/f_h] = \trans R\hlz{f_ht}$.
    If $f_0$ is $A$-regular
      and $f_1$,~\dots, $f_h$ is a d-sequence
      on~$A/ f_0^lA$ for all~$l > 0$,
      then $f_h$, $f_1/ f_h$,~\dots, $f_{h-1}/f_h$, $f_0$
      is a regular sequence on~$B$.
  \end{thm}

  \begin{pf}
    First we note that
      $f_1$,~\dots, $f_h$ is a d-sequence on~$A$.
    In fact,
      by using Krull's intersection theorem,
      we obtain
        \begin{align*}
          (f_1, \dots, f_{i-1}) \qtn f_i f_j
          & = \bigcap_{l =1}^\infty
            (f_0^l, f_1, \dots, f_{i-1}) \qtn f_i f_j
        \\
          & = \bigcap_{l =1}^\infty
            (f_0^l, f_1, \dots, f_{i-1}) \qtn f_j
        \\
          & =
            (f_1, \dots, f_{i-1}) \qtn f_j
        \end{align*}
      for any $1 \leq i \leq j \leq h$.
    Next we show that
        \begin{equation} \label{eqn:3.3.1}
          (f_1, \dots, f_{k-1}) \qtn f_k \cap \trans \idealq^n =
          (f_1, \dots, f_{k-1}) \trans \idealq^{n-1},
        \end{equation}
      for any $1 \leq k \leq h+1$ and $n>1$,
      where $f_{h+1} = 1$.
    If $a$ is an element of the left hand side,
      then $f_0^l a \in \idealq^n$
      for a sufficiently large~$l$.
    By~\cite[Theorem 1.3]{Goto-Yamagishi::theory},
      we have
        \begin{align*}
          f_0^l a &
          \in (f_1, \dots, f_{k-1}) \qtn f_k \cap \idealq^n
        \\
          &
          = (f_1, \dots, f_{k-1}) \idealq^{n-1}.
        \end{align*}
    Lemma~\ref{lem:3.2} says
          $$
            a \in [(f_1, \dots, f_{k-1}) \idealq^{n-1}] \qtn \angled{f_0}
            = (f_1, \dots, f_{k-1}) \trans \idealq^{n-1}.
          $$
    The inverse inclusion is clear.
    By~\eqref{eqn:3.3.1} and \cite[Theorem 1.7]{Goto-Yamagishi::theory},
      we obtain that
        $$
          f_h,
          \frac{f_1}{f_h},
          \dots,
          \frac{f_{h-1}}{f_h}
        $$
      is a regular sequence on~$B$.

    Finally we shall show that
      $f_0$ is regular on
      $B/(f_h, f_1/f_h, \dots, f_{h-1}/f_h)B$.
    Let $\alpha \in (f_h, f_1/f_h, \dots, f_{h-1}/f_h)B \qtn f_0$.
    For a sufficiently large~$n>1$,
      we may assume $\alpha = a_0/f_h^n$ and
        $$
          f_0 \frac{a_0}{f_h^n} =
          f_h \frac{a_h}{f_h^n} +
          \frac{f_1}{f_h} \frac{a_1}{f_h^n} +
          \dots +
          \frac{f_{h-1}}{f_h} \frac{a_{h-1}}{f_h^n}
        $$
      where $a_0$,~\dots, $a_h \in \trans\idealq^n$.
    Therefore
        $$
          f_h^{m+1} f_0 a_0 =
          f_h^m(f_h^2 a_h + f_1 a_1 + \dots + f_{h-1} a_{h-1})
        $$
      in~$A$ for some~$m>0$.
    Take an integer~$l$
      such that $f_0^l a_h \in \idealq^n$.
    Then
        \begin{align*}
          f_h^{m+2} f_0^l a_h &
          \in (f_0^{l+1}, f_1, \dots, f_{h-1}) \cap \idealq^{n+m+2}
        \\
          & = (f_0^{l+1}) \cap \idealq^{n+m+2}
          + (f_1, \dots, f_{h-1}) \idealq^{n+m+1}
        \\
          & \subseteq f_0^{l+1} \trans\idealq^{n+m+2}
          + (f_1, \dots, f_{h-1}) \idealq^{n+m+1}.
        \end{align*}
    If we put
        $$
          f_h^{m+2} f_0^l a_h =
          f_0^{l+1} b_0 + f_1 b_1 + \dots + f_{h-1} b_{h-1}
        $$
      where $b_0 \in \trans\idealq^{n+m+2}$ and
      $b_1$,~\dots, $b_{h-1} \in \idealq^{n+m+1}$,
      then
        \begin{align*}
          f_h^{m+2} a_h - f_0 b_0 & \in
          [(f_1, \dots, f_{h-1}) \idealq^{n+m+1}] \qtn \angled{f_0}
        \\
          & = (f_1, \dots, f_{h-1}) \trans\idealq^{n+m+1}.
        \end{align*}
    Let
        $$
          f_h^{m+2} a_h - f_0 b_0 =
          f_1 c_1 + \dots + f_{h-1} c_{h-1}
        $$
      where $c_1$,~\dots, $c_{h-1} \in \trans\idealq^{n+m+1}$.
    Then
        $$
          f_0(f_h^{m+1} a_0 - b_0) \in
          (f_1, \dots, f_{h-1}) \idealq^{n+m}.
        $$
    Therefore
        $$
          f_h^{m+1} a_0 - b_0 \in
          (f_1, \dots, f_{h-1}) \trans\idealq^{n+m},
        $$
      that is,
        $$
          \alpha - f_h \frac{b_0}{f_h^{n+m+2}}
          \in \left(
              \frac{f_1}{f_h},
              \dots,
              \frac{f_{h-1}}{f_h}
              \right)B.
        $$
    The proof is completed.
  \end{pf}

  In the rest of this section,
    we assume
    that $f_0$ is $A$-regular
    and that $f_1$,~\dots, $f_h$ is a u.s.d-sequence
    on~$A/ f_0^l A$ for all~$l>0$.
  Let $G = \bigoplus_{n \geq 0} \idealq^n/ \idealq^{n+1}$ and
      $
        \trans G = \bigoplus_{n \geq 0} \trans\idealq^n/
        \trans\idealq^{n+1}
      $
    be associated graded rings
    with respect to~$\idealq$
    and $\trans\idealq$,
    respectively.
  We shall compute local cohomology modules
    of~$\trans G$ and~$\trans R$
    with respect to~$\idealN = (f_0, \dots, f_h)R + R_+$.

  \begin{thm}
    If $p < h+1$, then
        $$
          [H_\idealN^p(\trans G)]_n = 0
          \quad
          \text{for $n \ne 1-p$}.
        $$
    Furthermore
        $$
          [H_\idealN^{h+1}(\trans G)]_n = 0
          \quad
          \text{for $n > -h$}.
        $$
  \end{thm}

  \begin{pf}
    We shall prove that
        \begin{equation} \label{eqn:3.4.1}
          [H_{(f_0, f_1t, \dots, f_kt)}^p(\trans G)]_n = 0
          \quad
          \text{for $n \ne 1-p$}
        \end{equation}
      if $p < k+1$
      by induction on~$k$.
    It is obvious that $f_0$ is $\trans G$-regular.
    Therefore $H_{(f_0)}^0(\trans G) = 0$.

    Suppose $k > 0$.
    Then $H_{(f_0, f_1t, \dots, f_{k-1}t)}^p(\trans G)_{f_kt} = 0$
      for $p < k$ by induction hypothesis.
    By Lemma~\ref{lem:2.2},
      we obtain isomorphisms
        $$
          H_{(f_0, f_1t, \dots, f_kt)}^p(\trans G) \cong
          H_{(f_0, f_1t, \dots, f_{k-1}t)}^p (\trans G)
          \quad
          \text{for $p<k$}.
        $$
    Therefore \eqref{eqn:3.4.1} is proved if $p < k$.
    We also obtain an exact sequence
        $$
          0 @>>>
          H_{(f_0, f_1t, \dots, f_kt)}^k(\trans G) @>>>
          H_{(f_0, f_1t, \dots, f_{k-1}t)}^k(\trans G) @>>>
          H_{(f_0, f_1t, \dots, f_{k-1}t)}^k(\trans G)_{f_kt}
        $$
      from Lemma~\ref{lem:2.2}.
    Hence $H_{(f_0, f_1t, \dots, f_kt)}^k(\trans G)$ is
      the limit of the direct system~$\{K_m\}_{m>0}$
      such that
        $$
          K_m =
          \frac
          {(f_0^m, (f_1t)^m, \dots, (f_{k-1}t)^m) \trans G
            \qtn \angled{f_kt}}
          {(f_0^m, (f_1t)^m, \dots, (f_{k-1}t)^m) \trans G}
          \, (m(k-1))
          \quad
          \text{for $m > 0$}
        $$
      and the homomorphism~$K_m \to K_{m'}$ is induced
      from the multiplication of
      $(f_0 \cdot f_1t \cdots f_{k-1}t)^{m' - m}$
      for any $m' > m$.
    We shall show that it is the zero map
      except for degree~$1-k$
      if $m'$ is sufficiently larger than~$m$.

    Let $\alpha$ be a homogeneous element of~$K_m$
      of degree~$n$
      and $a$ its representative.
    That is,
      $a \in \trans\idealq^{n+m(k-1)}$ and
        $$
          f_k^l a \in
          f_0^m \trans\idealq^{n+m(k-1)+l}
          + (f_1^m, \dots, f_{k-1}^m)
            \trans\idealq^{n+m(k-2)+l}
          + \trans\idealq^{n+m(k-1)+l+1}
        $$
      for some~$l>0$.
    Take an integer~$m' > m$
      such that $f_0^{m' - m}\trans \idealq \subseteq \idealq$.
    Then $f_0^{m' - m}\trans \idealq^n \subseteq \idealq^n$
      for any~$n>0$ by~\eqref{eqn:3.2.1}.
    By replacing $\alpha$ by its image in~$K_{m'}$,
      we may assume that $a \in \idealq^{n+m(k-1)}$ and
        $$
          f_k^l a \in
          f_0^m \trans\idealq^{n+m(k-1)+l}
          + (f_1^m, \dots, f_{k-1}^m) \idealq^{n+m(k-2) + l}
          + \idealq^{n+m(k-1) + l+1}.
        $$
    We put $f_k^l a = b + c$
      where
        $
          b \in f_0^m \trans\idealq^{n+m(k-1)+l}
          + (f_1^m, \dots, f_{k-1}^m) \idealq^{n+m(k-2) + l}
        $
      and $c \in \idealq^{n+m(k-1)+l+1}$.
    Then, by using~\cite[Theorem 2.6]{Goto-Yamagishi::theory},
      we obtain
        \begin{align*}
          c & \in (f_0^m, \dots, f_{k-1}^m, f_k^l)
            \cap \idealq^{n+m(k-1) + l+1}
        \\
          & \subseteq f_0^m \trans\idealq^{n+m(k-1) + l+1}
          + (f_1^m, \dots, f_{k-1}^m) \idealq^{n+m(k-2) + l+1}
          + f_k^l \idealq^{n+m(k-1) + 1}.
        \end{align*}
    If we put $c = b' + f_k^l a'$
      where
        $
          b' \in f_0^m \trans\idealq^{n+m(k-1) + l+1}
          + (f_1^m, \dots, f_{k-1}^m) \idealq^{n+m(k-2) + l+1}
        $
      and $a' \in \idealq^{n+m(k-1)+1}$,
      then $a-a'$ is also a representative of~$\alpha$.
    Therefore we may assume that $c=0$.

    By using~\cite[Theorem 2.8]{Goto-Yamagishi::theory},
      we obtain
        \begin{align*}
          a & \in (f_0^m, \dots, f_{k-1}^m) \qtn f_k \cap
            \idealq^{n+m(k-1)}
        \\
          & = (f_0^m) \cap \idealq^{n+m(k-1)}
          + (f_1^m, \dots, f_{k-1}^m) \idealq^{n+m(k-2)}
        \\
          & \quad +
            \sum
              \begin{Sb}
                I \subseteq \{1, \dots, k-1\}
              \\
                \sharp I \cdot (m-1) \geq n + m(k-1)
              \end{Sb}
            \left\{
              \prod_{i \in I}
                f_i^{m-1}
            \right\}
              \{
                [(f_0^m) + (f_i \mid i \in I)] \qtn f_k
              \}
        \\
          & \subseteq
          f_0^m \trans\idealq^{n+m(k-1)}
          + (f_1^m, \dots, f_{k-1}^m) \idealq^{n+m(k-2)}
          + \idealq^{n+m(k-1) + 1}
        \\
          & \quad +
            \sum
              \begin{Sb}
                I \subseteq \{1, \dots, k-1\}
              \\
                \sharp I \cdot (m-1) = n + m(k-1)
              \end{Sb}
            \left\{
              \prod_{i \in I}
                f_i^{m-1}
            \right\}
              \{
                [(f_0^m) + (f_i \mid i \in I)] \qtn f_k
              \}
        \end{align*}
    Here $\sharp I$ denotes the number of elements in~$I$.
    If $n > 1-k$,
      then there is no subset~$I$ of~$\{1, \dots, k-1\}$
      such that $\sharp I \cdot (m-1) = n+m(k-1)$.
    If $n < 1-k$,
      then such $I$ is a proper subset.
    Let $j \in \{1, \dots, k-1\} \setminus I$ and
        $$
          d \in [(f_0^m) + (f_i \mid i \in I)] \qtn f_k
          = [(f_0^m) + (f_i \mid i \in I)] \qtn f_j.
        $$
    Then
        $$
          (f_0 \cdots f_{k-1})
          \left\{ \prod_{i \in I} f_i^{m-1} \right\}
          d
          \in
          f_0^{m+1} \trans\idealq^{n+(m+1)(k-1)}
          + (f_1^{m+1}, \dots, f_{k-1}^{m+1}) \idealq^{n+(m+1)(k-2)}.
        $$
    In fact,
      if we put $f_j d = f_0^m e + g$
      where $g \in (f_i \mid i \in I)$,
      then $e \in \trans\idealq$.
    Thus the image of~$\alpha$ in~$K_{m+1}$ is zero
      if $n \ne 1-k$.

    Put $k=h$.
    Then
        $$
          [H_\idealN^p(\trans G)]_n =
          [H_{(f_0, f_1t, \dots, f_ht)}^p(\trans G)]_n =
          0
            \quad
          \text{for $n \ne 1 - p$}
        $$
      if $p < h+1$.
    The first assertion is proved.

    Next we compute $H_{(f_0, f_1t, \dots, f_ht)}^{h+1}(\trans G)$.
    It is the limit of the direct system~$\{K'_m\}_{m>0}$
      such that
        $$
          K'_m =
          \trans G/ (f_0^m, (f_1t)^m, \dots, (f_ht)^m) \trans G \, (mh)
          \quad
          \text{for $m>0$}
        $$
      and the homomorphism $K'_m \to K'_{m'}$ is induced
      from the multiplication of
      $(f_0 \cdot f_1t \cdots f_ht)^{m' - m}$
      for any $m' > m$.
    We shall show that it is the zero map for degree $n > -h$
      if $m'$ is sufficiently larger than~$m$.

    Let $\alpha$ be a homogeneous element of~$K'_m$ of degree~$n$
      and $a$ its representative.
    That is, $a \in \trans\idealq^{n+mh}$.
    If $n > -h$, then
        $$
          (f_0 \cdots f_h)^{m' - m} a \in
          \idealq^{n+m'h} =
          (f_1^{m'}, \dots, f_h^{m'}) \idealq^{n+m'(h-1)}
        $$
      for a sufficiently larger~$m'$ than~$m$.
    Thus the image of~$\alpha$ in~$K'_{m'}$ is zero
      if $n> -h$.
    Therefore $[H_\idealN^{h+1}(\trans G)]_n = 0$ for $n > -h$.
  \end{pf}

  By this theorem,
    we can compute local cohomology of~$\trans R$.

  \begin{cor} \label{cor:3.5}
    If $h=1$, $2$,
      then
        $$
          H_\idealN^p(\trans R) = 0
          \quad
          \text{for $p\ne 1$, $h+2$}
        $$
      and
        $
          H_\idealN^1(\trans R) =
          [H_\idealN^1(\trans R)]_0 =
          H_{(f_0, \dots, f_h)}^1(A)
        $.

    If $h \geq 3$,
      then
        $$
          H_\idealN^p(\trans R) = 0
          \quad
          \text{for $p=0$, $2$, $3$}
        $$
      and
        $
          H_\idealN^1(\trans R) =
          [H_\idealN^1(\trans R)]_0 =
          H_{(f_0, \dots, f_h)}^1(A)
        $.
    Furthermore,
       if $ 4 \leq p \leq h+1$,
      then
        $$
          [H_\idealN^p(\trans R)]_n =
          \begin{cases}
            H_{(f_0, \dots, f_h)}^{p-1}(A),
              & \text{for $-1 \geq n \geq 3-p$};
          \\
            0, & \text{otherwise}.
          \end{cases}
        $$
  \end{cor}

  \begin{pf}
    Passing through the completion,
      we may assume that $A$ possesses a dualizing complex.
    Since $H_\idealN^p(\trans G)$ is finitely graded for~$p < h+1$,
      $H_\idealN^p(\trans R)$ is finitely graded for~$p \leq h+1$
      \cite[Proposition~3]{Marley:94:finitely}.
    Considering the following two exact sequences
        $$
          0 @>>> \trans R_+ @>>> \trans R @>>> A @>>> 0
            \quad
            \text{and}
            \quad
          0 @>>> \trans R_+(1) @>>> \trans R @>>> \trans G @>>> 0,
        $$
      we obtain the assertion:
      see the proof of~\cite[Theorem 4.1]{Brodmann:84:local}.
  \end{pf}

  Let $S = \trans R/ R$, that is,
    $S = \bigoplus_{n>0} \trans\idealq^n / \idealq^n$.
  The following proposition shall play an important role
    in the next section.

  \begin{prop} \label{prop:3.6}
    If $p < h$, then
        $$
          [H_\idealN^p(S)]_n =0
          \quad
          \text{for $n \ne 1-p$.}
        $$
    Moreover,
        $$
          [H_\idealN^h(S)]_n = 0
          \quad
          \text{for $n > 1-h$}.
        $$
  \end{prop}

  \begin{pf}
    In the same way as the proof of Theorem~\ref{thm:3.3},
      we find that $f_1$,~\dots, $f_h$ is a u.s.d-sequence on~$A$.
    Hence, by using \cite[Theorem 4.2]{Goto-Yamagishi::theory},
        $$
          [H_{(f_1t, \dots, f_ht)}^p(G)]_n = 0
          \quad
          \text{for $n \ne -p$}
        $$
      if $p<h$.
    Furthermore,
        $$
          [H_{(f_1t, \dots, f_ht)}^h(G)]_n = 0
          \quad
          \text{for $n> -h$}.
        $$
    By using Lemma~\ref{lem:2.2},
      we obtain
        $$
          [H_\idealN^p(G)]_n = 0
          \quad
          \text{for $n \ne 1-p$, $-p$}
        $$
      if $p < h$ and
        $$
          [H_\idealN^p(G)]_n = 0
          \quad
          \text{for $n>1-p$}
        $$
      if $p = h$, $h+1$.

    Since
        $
          \trans\idealq^2 =
          \idealq \trans\idealq
        $,
      there exists an exact sequence
        $$
          0 @>>> S(1) @>>> G @>\phi>> \trans G @>>> S @>>> 0.
        $$
    Let $T$ be the image of~$\phi$.
    We shall show
        $$
          [H_\idealN^p(S)]_n =
          [H_\idealN^p(T)]_n = 0
          \quad
          \text{for $n> 1-p$}
        $$
      by induction on~$h-p$.
    If $p > h+1$, then the assertion is obvious.
    Let $p \leq h+1$.
    Then following two exact sequences
        \begin{gather*}
          H_\idealN^p(\trans G) @>>>
          H_\idealN^p(S) @>>>
          H_\idealN^{p+1}(T) @>>>
          H_\idealN^{p+1}(\trans G),
        \\
          H_\idealN^p(G) @>>>
          H_\idealN^p(T) @>>>
          H_\idealN^{p+1}(S)(1) @>>>
          H_\idealN^{p+1}(G)
        \end{gather*}
      and the induction hypothesis imply
        $$
          [H_\idealN^p(S)]_n =
          [H_\idealN^p(T)]_n = 0
          \quad
          \text{for $n>1-p$}.
        $$

    In the same way,
      we can prove that
        $$
          [H_\idealN^p(S)]_n =
          [H_\idealN^p(T)]_n = 0
          \quad
          \text{for $n< 1-p$}
        $$
      if $p< h$
      by induction on~$p$.
  \end{pf}

  Finally we show that
    $\trans R$ is an ideal transform of~$R$
    in a sense.

  \begin{prop}
    $\trans R_+ = D_{(f_0, \dots, f_h)}^0(R_+)$.
  \end{prop}

  \begin{pf}
    We first show that
      $f_0$, $f_1$ is a regular sequence on~$\trans R_+$.
    Let $n>0$.
    Since $f_0$ is $A$-regular,
      it is also $\trans\idealq^n$-regular.
    Let
        $
          a \in
          [f_0 \trans\idealq^n] \qtn f_1 \cap
          \trans \idealq^n
        $.
    Then $f_0^l a \in \idealq^n$ for a sufficiently large~$l$.
    Since $f_1 a \in (f_0)$,
      we have
        $
          f_0^l a \in
          (f_0^{l+1}) \qtn f_1 \cap \idealq^n \subseteq
          f_0^{l+1} \trans\idealq^n
        $,
      that is, $a \in f_0 \trans\idealq^n$.
    Thus we have shown that
      $f_1$ is $\trans R_+ / f_0 \trans R_+$-regular.

    By this and \eqref{eqn:2.1.1},
      we obtain
        \begin{equation} \label{eqn:3.7.1}
          D_{(f_0, \dots, f_h)}^0(R_+) \subseteq
          D_{(f_0, \dots, f_h)}^0(\trans R_+) =
          \trans R_+.
        \end{equation}
    Since
        $
          \trans\idealq^n =
          \idealq^{n-1} \trans\idealq
        $
      for $n \geq 2$,
      $(f_0^l, f_1, \dots, f_h) \trans R_+ \subseteq R_+$
      for a sufficiently large~$l$.
    Hence, we obtain the inverse inclusion of~\eqref{eqn:3.7.1}.
    The proof is completed.
  \end{pf}

\setcounter{equation}{0}
\section{%
  A blowing-up with respect to a certain subsystem of parameters}
    \label{sec:4}

  In this section,
    we assume that $A$ possesses a dualizing complex.
  We fix an integer $s \geq \dim A/ \ideala_A(A)$.
  Since $\dim A/ \ideala_A(M) < \dim M$
    for any finitely generated $A$-module~$M$
    \cite[Korollar 2.2.4]{Schenzel:82:dualisierende},
    there exists a system of parameters
    $x_1$,~\dots, $x_d$ for~$A$ such that
      \begin{equation} \label{eqn:4.0.1}
        \begin{cases}
          x_{s+1}, \dots, x_d \in \ideala_A(A);
        \\
          x_i \in \ideala_A(A/(x_{i+1}, \dots, x_d)),
            & \text{for $i \leq s$}.
        \end{cases}
      \end{equation}
    This notion is a slight improvement of
      a p-standard system of parameters,
      which was introduced by Cuong~\cite{Cuong:91:dimension}.
    He also gave the statement~(1) of Theorem~\ref{thm:4.2}.
    The author was personally taught it by him.

  \begin{lem} \label{lem:4.1}
    Let $n_1$,~\dots, $n_i$ be arbitrary positive integers.
    Then
        \begin{multline*}
          (x_1^{n_1}, \dots, x_{i-1}^{n_{i-1}},
            x_{k+1}, \dots, x_d) \qtn x_i^{n_i} \cap
          (x_1^{n_1}, \dots, x_{i-1}^{n_{i-1}}, x_k, \dots, x_d)
        \\
          = (x_1^{n_1}, \dots, x_{i-1}^{n_{i-1}}, x_{k+1}, \dots, x_d)
        \end{multline*}
      for any $1 \leq i \leq k \leq d$.
  \end{lem}

  \begin{pf}
    It is obvious that the left hand side
      contains the right one.
    Let $a$ be an element of the left hand side
      and $a = b + x_k c$ where
        $
          b \in (x_1^{n_1}, \dots, x_{i-1}^{n_{i-1}},
          x_{k+1}, \dots, x_d)
        $.
    Then
        \begin{align*}
          c & \in (x_1^{n_1}, \dots, x_{i-1}^{n_{i-1}},
            x_{k+1}, \dots, x_d) \qtn x_i^{n_i} x_k
        \\
          & = (x_1^{n_1}, \dots, x_{i-1}^{n_{i-1}},
            x_{k+1}, \dots, x_d) \qtn x_k
        \end{align*}
      by Lemma~\ref{lem:2.4}.
    Therefore
        $
          x_k c, a \in (x_1^{n_1}, \dots, x_{i-1}^{n_{i-1}},
            x_{k+1}, \dots, x_d)
        $.
    The proof is completed.
  \end{pf}

  Let $\idealq = (x_{s+1}, \dots, x_d)$.
  Lemma~\ref{lem:2.2} assures us
    that $x_{s+1}$, \dots, $x_d$ is a u.s.d-sequence on~$A$.
  Furthermore, we have the following theorem:

  \begin{thm} \label{thm:4.2}
    \rom{(1)}
      The sequences $x_1^{n_1}$,~\dots, $x_s^{n_s}$,
        $x_{\sigma(s+1)}^{n_{s+1}}$,~\dots, $x_{\sigma(d)}^{n_d}$ is
        a d-sequence on~$A$
        for any positive integers $n_1$,~\dots, $n_d$
        and for any permutation~$\sigma$ on $s+1$,~\dots, $d$.

    \rom{(2)}
      If $s>0$, then $x_1^{n_1}$,~\dots, $x_s^{n_s}$ is
        a d-sequence on~$A/ \idealq^n$
        for any positive integers $n_1$,~\dots, $n_s$
        and~$n$.
  \end{thm}

  \begin{pf}
    (1):~Let $1 \leq i \leq j \leq d$.
    We have only to prove that
        $$
          (x_1^{n_1}, \dots, x_{i-1}^{n_{i-1}}) \qtn x_i^{n_i} x_j^{n_j}
          =
          (x_1^{n_1}, \dots, x_{i-1}^{n_{i-1}}) \qtn x_j^{n_j}
        $$
      for any positive integers $n_1$,~\dots, $n_d$.
    If $j>s$, then the both sides are equal to
      $(x_1^{n_1}, \dots, x_{i-1}^{n_{i-1}}) \qtn \ideala_A(A)$.

    Assume that $j \leq s$
      and take an element~$a$ of the left hand side.
    By using Lemma~\ref{lem:2.4}, we get
        \begin{align*}
          a & \in (x_1^{n_1}, \dots, x_{i-1}^{n_{i-1}},
            x_{j+1}, \dots, x_d) \qtn x_i^{n_i} x_j^{n_j}
        \\
          & = (x_1^{n_1}, \dots, x_{i-1}^{n_{i-1}},
            x_{j+1}, \dots, x_d) \qtn x_j^{n_j}.
        \end{align*}
      Hence we have
          \begin{align*}
            x_j^{n_j} a & \in (x_1^{n_1}, \dots, x_{i-1}^{n_{i-1}})
                \qtn x_i^{n_i} \cap
            (x_1^{n_1}, \dots, x_{i-1}^{n_{i-1}},
                x_{j+1}, \dots, x_d)
          \\
            & = (x_1^{n_1}, \dots, x_{i-1}^{n_{i-1}})
          \end{align*}
        by repeating to use Lemma~\ref{lem:4.1}.

    (2):~If $n=1$,
      then the assertion is proved in the same way as above.
    Let $1 \leq i \leq j \leq s$ and $n>1$.
    Then $x_{s+1}$,~\dots, $x_d$ is a d-sequence
      on~$A/(x_1^{n_1}, \dots, x_{i-1}^{n_{i-1}}, x_i^{n_i} x_j^{n_j})$.
    By using Lemma~\ref{lem:3.2},
      we obtain
        \begin{align*}
          \lefteqn{[(x_1^{n_1}, \dots, x_{i-1}^{n_{i-1}}) + \idealq^n]
            \qtn x_i^{n_i} x_j^{n_j}}
            \qquad
        \\
          & = (x_1^{n_1}, \dots, x_{i-1}^{n_{i-1}})
            \qtn x_i^{n_i} x_j^{n_j}
          + \idealq^{n-1} [(x_1^{n_1}, \dots, x_{i-1}^{n_{i-1}},
            x_{s+1}, \dots, x_d) \qtn x_i^{n_i} x_j^{n_j}]
        \\
          & = (x_1^{n_1}, \dots, x_{i-1}^{n_{i-1}})
            \qtn x_j^{n_j}
          + \idealq^{n-1} [(x_1^{n_1}, \dots, x_{i-1}^{n_{i-1}},
            x_{s+1}, \dots, x_d) \qtn x_j^{n_j}]
        \\
          & \subseteq
            [(x_1^{n_1}, \dots, x_{i-1}^{n_{i-1}}) + \idealq^n]
            \qtn x_j^{n_j}.
        \end{align*}
    Here the second equality follows from the case of~$n=1$.
    Thus the proof is completed.
  \end{pf}

  In the same way as the proof of Theorem~\ref{thm:3.3},
    we find that any subsequence of $x_1^{n_1}$,~\dots, $x_d^{n_d}$
    is a d-sequence on~$A$
    and any subsequence of $x_1^{n_1}$,~\dots, $x_s^{n_s}$
    is a d-sequence on~$A/ \idealq^n$
    for arbitrary positive integers $n_1$,~\dots, $n_d$ and~$n$.

  \begin{cor}
    Fix an integer~$k$ such that $1 \leq k \leq d$.
    Then
        $$
          H_{(x_k, \dots, x_d)}^p(A) =
          \dlim_m \frac{(x_k^m, \dots, x_{k+p-1}^m) \qtn x_{k+p}}
            {(x_k^m, \dots, x_{k+p-1}^m)}
          \quad
          \text{for $p< d-k+1$}.
        $$
  \end{cor}

  \begin{pf}
    We shall prove that
        $$
          H_{(x_k, \dots, x_l)}^p(A) =
          \dlim_m \frac{(x_k^m, \dots, x_{k+p-1}^m) \qtn x_{k+p}}
            {(x_k^m, \dots, x_{k+p-1}^m)}
          \quad
          \text{for $p< l-k+1$}
        $$
      by induction on~$l \geq k$.
    If $l=k$, then $H_{(x_k)}^0(A) = 0 \qtn_A x_k$.

    Suppose $l > k$.
    Then $x_k$,~\dots, $x_{l-1}$ is a regular sequence on~$A_{x_l}$
      because $x_k$,~\dots, $x_l$ is a d-sequence on~$A$.
    Hence we obtain isomorphisms
        $$
          H_{(x_k, \dots, x_l)}^p(A) \cong
          H_{(x_k, \dots, x_{l-1})}^p(A)
          \quad
          \text{for all $p<l-k$}
        $$
      and an exact sequence
        $$
          0 @>>>
          H_{(x_k, \dots, x_l)}^{l-k}(A) @>>>
          H_{(x_k, \dots, x_{l-1})}^{l-k}(A) @>>>
          H_{(x_k, \dots, x_{l-1})}^{l-k}(A)_{x_l}
        $$
      by Lemma~\ref{lem:2.2}.
    This exact sequence
      is the direct limit of the exact sequence
        $$
          0 @>>>
          \frac{(x_k^m, \dots, x_{l-1}^m) \qtn x_l}
            {(x_k^m, \dots, x_{l-1}^m)} @>>>
          A/(x_k^m, \dots, x_{l-1}^m) @>>>
          [A/(x_k^m, \dots, x_{l-1}^m)]_{x_l}
        $$
    Thus the proof is completed.
  \end{pf}

  If $s=0$,
    then $\proj A[\idealq t] \rightarrow \spec A$ is
    a Macaulayfication of~$\spec A$:
    see Theorem~\ref{thm:5.1} for details.
  In the rest of this section,
    we shall observe $\proj A[\idealq t]$
    when $s>0$.
  Assume that $s>0$ and fix an integer~$k$
    such that $1 \leq k \leq s$.
  We shall compute local cohomology modules of $R = A[\idealq t]$
    with respect to~$(x_k, \dots, x_{s+1})$.
  Let $\idealM = \idealm R + R_+$.

  \begin{thm}
    $H_{(x_k, \dots, x_{s+1})}^0(R) = 0 \qtn_A x_k$.
  \end{thm}

  \begin{pf}
    Since $x_k$, $x_{s+1}$,~\dots, $x_d$ is a d-sequence on~$A$,
      $0 \qtn_A x_k \cap \idealq^n = 0$ for~$n>0$
      by~\cite[Theorem 1.3]{Goto-Yamagishi::theory}.
    That is,
        $$
          H_{(x_k, \dots, x_{s+1})}^0(\idealq^n) =
          \begin{cases}
            0 \qtn_A x_k, & \text{if $n=0$};
          \\
            0, & \text{otherwise}.
          \end{cases}
        $$
    Therefore,
        $
          H_{(x_k, \dots, x_{s+1})}^0(R) =
          \bigoplus_{n \geq 0} H_{(x_k, \dots, x_{s+1})}^0(\idealq^n)
          = 0 \qtn_A x_k
        $.
  \end{pf}

  Let $C = A[t]/R$,
    that is,
    $C = \bigoplus_{n>0} A/ \idealq^n$.

  \begin{lem}
    For $k \leq l \leq s+1$ and $p \leq l-k$,
      the natural homomorphism
        $$
          \alpha_l^p \colon
          H_{(x_k, \dots, x_l)}^p(A[t]) @>>>
          H_{(x_k, \dots, x_l)}^p(C)
        $$
      is a monomorphism except for degree~$0$.
  \end{lem}

  \begin{pf}
    We shall work by induction on~$l$.
    If $l=k$,
      then $0 \qtn_A x_k \cap \idealq^n = 0$
      for $n>0$.
    Therefore
        $$
          \alpha_k^0 \colon
          0 \qtn_{A[t]} x_k @>>>
          \bigoplus_{n>0} \idealq^n \qtn x_k / \idealq^n
        $$
      is a monomorphism except for degree~$0$.
    Let $k < l \leq s$.
    Then $x_k$,~\dots, $x_{l-1}$ is a regular sequence
      on~$A_{x_l}$ and on~$C_{x_l}$ by Theorem~\ref{thm:4.2}.
    By using Lemma~\ref{lem:2.2},
      we obtain commutative diagrams
        $$
          \begin{CD}
            H_{(x_k, \dots, x_l)}^p(A[t])
            @>\sim>>
            H_{(x_k, \dots, x_{l-1})}^p(A[t])
          \\
            @V{\alpha_l^p}VV @V{\alpha_{l-1}^p}VV
          \\
            H_{(x_k, \dots, x_l)}^p(C)
            @>\sim>>
            H_{(x_k, \dots, x_{l-1})}^p(C)
          \end{CD}
          \quad
          \text{for $p<l-k$}
        $$
      and
        $$
          \begin{CD}
            0 @>>>
            H_{(x_k, \dots, x_l)}^{l-k} (A[t]) @>>>
            H_{(x_k, \dots, x_{l-1})}^{l-k} (A[t]) @>>>
            H_{(x_k, \dots, x_{l-1})}^{l-k} (A[t])_{x_l}
          \\
            @. @V{\alpha_l^{l-k}}VV @VVV @VVV
          \\
            0 @>>>
            H_{(x_k, \dots, x_l)}^{l-k} (C) @>>>
            H_{(x_k, \dots, x_{l-1})}^{l-k} (C) @>>>
            H_{(x_k, \dots, x_{l-1})}^{l-k} (C)_{x_l}
          \end{CD}
        $$
      whose rows are exact.
    Therefore the assertion is true for $p<l-k$
      and we find that
      $\alpha_l^{l-k}$ is the direct limit of
        $$
          \alpha_{l,m} \colon
          \frac{(x_k^m, \dots, x_{l-1}^m) A[t] \qtn x_l}
            {(x_k^m, \dots, x_{l-1}^m) A[t]} @>>>
          \bigoplus_{n>0}
            \frac{[(x_k^m, \dots, x_{l-1}^m) + \idealq^n] \qtn x_l}
            {(x_k^m, \dots, x_{l-1}^m) + \idealq^n}.
        $$
    Since $x_l$, $x_{s+1}$,~\dots, $x_d$ is a d-sequence
      on~$A/(x_k^m, \dots, x_{l-1}^m)$,
        $$
          (x_k^m, \dots, x_{l-1}^m) \qtn x_l \cap
          [(x_k^m, \dots, x_{l-1}^m) + \idealq^n] =
          (x_k^m, \dots, x_{l-1}^m)
          \quad
          \text{for $n>0$}.
        $$
    Therefore $\alpha_{l,m}$ is a monomorphism except for degree~$0$
      and $\alpha_l^{l-k}$ is also.

    If $l=s+1$,
      then $x_k$,~\dots, $x_s$ is a regular sequence on~$A_{x_{s+1}}$
      and $C_{x_{s+1}} = 0$.
    The assertion is proved in the same way as above.
  \end{pf}

  Of course, $\alpha_{s+1}^p$ is the zero map
    in degree~$0$.
  Therefore there exists an exact sequence
      \begin{equation} \label{eqn:4.5.1}
        0 @>>>
        \coker \alpha_{s+1}^{p-1} @>>>
        H_{(x_k, \dots, x_{s+1})}^p(R) @>>>
        H_{(x_k, \dots, x_{s+1})}^p(A) @>>>
        0
      \end{equation}
    for $0 < p \leq s-k+1$.

  \begin{thm} \label{thm:4.6}
    Let $0 \leq q \leq s-k$.
    Then
        $$
          (x_{k+q}, \dots, x_d) \coker \alpha^q_{s+1} = 0
        $$
      and $H_\idealM^p(\coker \alpha^q_{s+1})$ is finitely graded
      for~$p< d-s$.
  \end{thm}

  \begin{pf}
    We know that
        $
          \coker \alpha_{s+1}^q =
          \coker \alpha_{k+q}^q =
          \dlim_m \coker \alpha_{k+q,m}
        $
      and
        $$
          \coker \alpha_{k+q,m} =
          \bigoplus_{n>0}
          \frac{[(x_k^m, \dots, x_{k+q-1}^m) + \idealq^n] \qtn x_{k+q}}
          {(x_k^m, \dots, x_{k+q-1}^m) \qtn x_{k+q} + \idealq^n}.
        $$
    By using~Theorem~\ref{thm:4.2} and Lemma~\ref{lem:3.2},
      we obtain
        \begin{multline} \label{eqn:4.6.1}
          [(x_k^m, \dots, x_{k+q-1}^m) + \idealq^n] \qtn x_{k+q}
        \\
          = (x_k^m, \dots, x_{k+q-1}^m) \qtn x_{k+q} +
          \idealq^{n-1} [(x_k^m, \dots, x_{k+q-1}^m,
            x_{s+1}, \dots, x_d) \qtn x_{k+q}].
        \end{multline}
    Therefore $\coker \alpha_{k+q,m}$ is annihilated
      by~$(x_{k+q}, \dots, x_d)$
      and $\coker \alpha_{s+1}^q$ is also.

    Next we compute local cohomology modules
      of~$\coker \alpha_{s+1}^q$.
    We note that $x_{k+q}$ is a regular element
      on~$A/(x_k^m, \dots, x_{k+q-1}^m) \qtn x_{k+q}$
      and that $x_{s+1}$,~\dots, $x_d$ is a u.s.d-sequence
      on~$A/ (x_k^m, \dots, x_{k+q-1}^m) \qtn x_{k+q} + (x_{k+q}^l)$
      for any $l>0$:
      see \cite[Proposition 2.2]{Huneke:82:theory}.
    Therefore, by Proposition~\ref{prop:3.6},
        \begin{equation} \label{eqn:4.6.2}
          \text{
              $
                H_{(x_{k+q}, x_{s+1}, \dots, x_d)R + R_+}^p
                  (\coker \alpha_{k+q, m})
              $
            is concentrated in degree~$1-p$}
        \end{equation}
      if $p<d-s$.
    Hence
        $
          H_{(x_{k+q}, x_{s+1}, \dots, x_d)R + R_+}^p
          (\coker \alpha_{s+1}^q)
        $
      is also.
    By the spectral sequence
        $
          E_2^{pq} =
          H_\idealM^p H_{(x_{k+q}, x_{s+1}, \dots, x_d)R + R_+}^q (-)
          \Rightarrow
          H_\idealM^{p+q}(-)
        $,
      we obtain the second assertion.
  \end{pf}

  Next we compute $H_{(x_k, \dots, x_{s+1})}^{s-k+2}(R)$.

  \begin{thm}
    Let $A_m = A/(x_k^m, \dots, x_s^m)$ and
      $\idealq_m = \idealq A_m$
      for any positive integer~$m$.
    Then
        $$
          H_{(x_k, \dots, x_{s+1})}^{s-k+2}(R) =
          \dlim_{m,l} A_m[\idealq_m t]/ x_{s+1}^l A_m[\idealq_m t].
        $$
    In particular,
      $H_\idealM^p H_{(x_k, \dots, x_{s+1})}^{s-k+2}(R)$ is
      finitely graded for~$p<d-s$.
  \end{thm}

  \begin{pf}
    We consider the exact sequence
      \begin{multline*}
        H_{(x_k, \dots, x_s)}^{s-k}(A[t]) @>\alpha_s^{s-k}>>
        H_{(x_k, \dots, x_s)}^{s-k}(C) @>>>
        H_{(x_k, \dots, x_s)}^{s-k+1}(R) @>>>
      \\
        @>>>
        H_{(x_k, \dots, x_s)}^{s-k+1}(A[t]) @>\beta>>
        H_{(x_k, \dots, x_s)}^{s-k+1}(C).
      \end{multline*}
    Since $\beta$ is the direct limit of
        $$
          A[t]/(x_k^m, \dots, x_s^m) A[t] @>>>
          C/ (x_k^m, \dots, x_s^m) C,
        $$
      we have $\ker \beta = \dlim_m A_m [\idealq_m t]$.
    Taking local cohomology modules of a short exact sequence
        $$
          0 @>>>
          \coker \alpha_s^{s-k} @>>>
          H_{(x_k, \dots, x_s)}^{s-k+1}(R) @>>>
          \ker \beta @>>>
          0
        $$
      with respect to $(x_{s+1})$,
      we obtain
        \begin{equation} \label{eqn:4.7.1}
          H_{(x_{s+1})}^1 H_{(x_k, \dots, x_s)}^{s-k+1}(R) =
          H_{(x_{s+1})}^1 (\ker \beta),
        \end{equation}
      because $\coker \alpha_s^{s-k} = \coker \alpha_{s+1}^{s-k}$
      is annihilated by~$x_{s+1}$.
    The left hand side of~\eqref{eqn:4.7.1} coincides
      with~$H_{(x_k, \dots, x_{s+1})}^{s-k+2}(R)$
      by Lemma~\ref{lem:2.2}.
    Thus the first assertion is proved.

    Since $x_{s+1}$,~\dots, $x_d$ is a u.s.d-sequence on~$A_m$,
      $H_{(x_{s+1}, \dots, x_d)R + R_+}^p(A_m [\idealq_m t])$ is
      concentrated in degree $0 \geq n \geq s-d+2$
      if $p \leq d-s$:
      see \cite[Theorem 4.1]{Goto-Yamagishi::theory}.
    From the exact sequence
        $$
          0 @>>>
          0 \qtn_{A_m} x_{s+1} @>>>
          A_m [\idealq_m t] @>x_{s+1}^l>>
          A_m [\idealq_m t] @>>>
          A_m[\idealq_m t] / x_{s+1}^l A_m [\idealq_m t] @>>> 0
        $$
      and the spectral sequence
        $
          E_2^{pq} =
          H_\idealM^p H_{(x_{s+1}, \dots, x_d)R + R_+}^q(-)
          \Rightarrow
          H_\idealM^{p+q}(-)
        $,
      we find that
        \begin{equation} \label{eqn:4.7.2}
          \text{
            $H_\idealM^p(A_m[\idealq_m t]/ x_{s+1}^l A_m [\idealq_m t])$
            is concentrated in degree $0 \geq n \geq s-d+2$}
        \end{equation}
      if $p< d-s$.
    Taking the direct limit of it,
      we obtain the second assertion.
  \end{pf}

  Finally we compute
    local cohomology modules of
    $B = A[\idealq/ x_{s+1}] = R\hlz{x_{s+1}t}$.

  \begin{thm} \label{thm:4.8}
    Let $\idealn$ be a maximal ideal of~$B$.
    Then
        $$
          H_\idealn^p H_{(x_k, \dots, x_{s+1})}^q(B) = 0
          \quad
          \text{if $q=0$ or $p<d-s-1$.}
        $$
    Furthermore
      $(x_{k+q-1}, \dots, x_{s+1}) H_{(x_k, \dots, x_{s+1})}^q(B) =0$
      for~$q < s-k+2$.
  \end{thm}

  \begin{pf}
    Since the blowing-up $\proj R \to \spec A$ is a closed map,
      there exists a homogeneous prime ideal~$\idealp$ of~$R$
      such that $x_{s+1}t \notin \idealp$,
      $\dim R/ \idealp = 1$ and
      $\idealn = [\idealp R_{x_{s+1}t}]_0$.

    Since $x_{s+1}$ is $B$-regular,
      $H_{(x_k, \dots, x_{s+1})}^0(B) =0$.

    Let $1 \leq q \leq s-k+1$.
    By applying Lemma~\ref{lem:2.5}
      to~\eqref{eqn:4.6.2},
      we obtain
        $$
          H_\idealn^p((\coker \alpha_{k+q-1,m})\hlz{x_{s+1}t}) =0
          \quad
          \text{for $p<d-s-1$}.
        $$
    By taking the direct limit of it
      and using~\eqref{eqn:4.5.1},
      we have
        $$
          H_\idealn^p H_{(x_k, \dots, x_{s+1})}^q(B) = 0
          \quad
          \text{for $p<d-s-1$}.
        $$
    Moreover Theorem~\ref{thm:4.6}
      also assures us
        $
          (x_{k+q-1}, \dots, x_{s+1}) H_{(x_k, \dots, x_{s+1})}^q(B) =0
        $.

    Next we consider $H_{(x_k, \dots, x_{s+1})}^{s-k+2}(B)$.
    By applying Lemma~\ref{lem:2.5} to~\eqref{eqn:4.7.2} and
      by taking direct limit,
      we have
        $$
          H_\idealn^p H_{(x_k, \dots, x_{s+1})}^{s-k+2}(B) = 0
          \quad
          \text{for $p<d-s-1$}.
        $$
    Thus the proof is completed.
  \end{pf}

\section{Macaulayfications of local rings}
  \label{sec:5}

  In this section,
    we shall construct a Macaulayfication of
    the affine scheme~$\spec A$
    if its non-Cohen-Macaulay locus is of dimension~$2$.
  Assume that $A$ possesses a dualizing complex
    and $\dim A/ \idealp = d$
    for any associated prime ideal~$\idealp$ of~$A$.
  Then $V(\ideala_A(A))$ coincides with
    the non-Cohen-Macaulay locus of~$A$.
  We fix an integer $s \geq \dim A/ \ideala_A(A)$
    and
    let $x_1$,~\dots, $x_d$ be a system of parameters for~$A$
    satisfying~\eqref{eqn:4.0.1}.

  First we refine Faltings' results~%
    \cite[Satz 2, 3]{Faltings:78:Macaulay}.
  Let $\idealq = (x_{s+1}, \dots, x_d)$,
    $R = A[\idealq t]$ and
    $X = \proj R$.

  \begin{thm} \label{thm:5.1}
    With notation as above,
        $$
          \depth \sheafO_{X,p} \geq d-s
          \quad
          \text{for any closed point~$p$ of~$X$.}
        $$
    If $s=0$ or $A/ \idealq$ is Cohen-Macaulay,
      then $X$ is Cohen-Macaulay.
  \end{thm}

  \begin{pf}
    Since $x_{s+1}$,~\dots, $x_d$ is a u.s.d-sequence on~$A$,
      $H_{(x_{s+1}, \dots, x_d)R + R_+}^p(R)$
      is finitely graded for $p \leq d-s$:
      see \cite[Theorem 4.1]{Goto-Yamagishi::theory}.
    By using Lemma~\ref{lem:2.5},
      we obtain the first assertion.

    Furthermore
      since $\dim \sheafO_{X,p} =d$ for any closed point~$p$ of~$X$,
      $X$ is Cohen-Macaulay if $s=0$.

    Assume that $s>0$ and $A / \idealq$ is Cohen-Macaulay.
    Then $x_1$,~\dots, $x_s$ is a regular sequence on~$A/ \idealq$.
    We use theorems in Section~\ref{sec:4}
      as $k=1$.
    From~\eqref{eqn:4.6.1},
      we find that $\coker \alpha_{s+1}^q =0$ for all~$q \leq s-1$.
    That is,
      $H_\idealM^p H_{(x_1, \dots, x_{s+1})}^q(R)$ is
      finitely graded
      if $p<d-s$ or $q<s+1$.
    By the spectral sequence~%
        $
          E_2^{pq} =
          H_\idealM^p H_{(x_1, \dots, x_{s+1})}^q(-) \Rightarrow
          H_\idealM^{p+q} (-)
        $,
      we find that $H_\idealM^p(R)$ is
      finitely graded for $p<d+1$.
    Lemma~\ref{lem:2.5} assures us
        $$
          \depth \sheafO_{X,p} \geq d
          \quad
          \text{for any closed point~$p$ of~$X$}.
        $$
    The proof is completed.
  \end{pf}

  From now on,
    we assume that $s > 0$.

  Since $x_s$ is $A$-regular,
    $\idealq$ is a reduction of~$\trans\idealq = \idealq \qtn x_s$
    by~\eqref{eqn:3.2.1}.
  We put $\trans R = A[\trans\idealq t]$ and
    $\trans X = \proj \trans R$.
  Then $\trans X \rightarrow X$ is a finite morphism.

  \begin{thm}
    With notation as above,
        $$
          \depth \sheafO_{\trans X, \trans p} \geq d-s+1
          \quad
          \text{for any closed point~$\trans p$ of~$\trans X$}
        $$
    In particular, if $s=1$,
      then $\trans X$ is Cohen-Macaulay.
  \end{thm}

  \begin{pf}
    By Corollary~\ref{cor:3.5},
      $H_{(x_s, \dots, x_d)R + R_+}^p(\trans R)$ is finitely graded
      for $p \leq d-s+1$.
    By using Lemma~\ref{lem:2.5},
      we obtain the assertion.
  \end{pf}

  Next we consider an ideal~%
      $
        \idealb =
        \idealq^2 + x_s \idealq = (x_s, \dots, x_d) \idealq
      $.
  We put $S = A[\idealb t]$ and $Y = \proj S$.
  Then $Y$ is the blowing-up of~$X$
    with center~$(x_s, \dots, x_d) \sheafO_X$.

  \begin{thm} \label{thm:5.3}
    With notation as above,
        $$
          \depth \sheafO_{Y,q} \geq d-s+1
          \quad
          \text{for any closed point~$q$ of~$Y$}.
        $$
    Furthermore,
      if $s=1$ or $A$ is Cohen-Macaulay,
      then $Y$ is Cohen-Macaulay.
  \end{thm}

  \begin{pf}
    Since
        $
          (x_s x_{s+1}, \dots, x_s x_d, x_{s+1}^2, \dots, x_d^2)
            \idealb^{d-s-1}
          = \idealb^{d-s}
        $,
      we have only to compute the depth of
      $C_0 = A[\idealb/ x_s x_{s+1}]$ and $C_1 = A[\idealb/ x_{s+1}^2]$.
    If we put $B = A[\idealq/ x_{s+1}]$,
      then
        \begin{align*}
          C_0 & = B[x_{s+1}/ x_s]
            \cong B[T]/(x_sT - x_{s+1}) \qtn \angled{x_s},
        \\
          C_1 & = B[x_s/ x_{s+1}]
            \cong B[T]/(x_{s+1}T - x_s) \qtn \angled{x_{s+1}},
        \end{align*}
      where $T$ denotes an indeterminate.
    We note that $B$, $C_0$, $C_1$ are
      subrings of the total quotient ring of~$A$
      because $x_1$,~\dots, $x_d$ are $A$-regular elements.

    First we consider $C_0$.
    We regard it as a homomorphic image of~$B[T]$.
    Let $\ideall_0$ be a maximal ideal of~$C_0$ and
      $\idealn = \ideall_0 \cap B$.
    Then $\idealn$ is a maximal ideal of~$B$
      because $\spec C_0 \cup \spec C_1 \rightarrow \spec B$ is
      a blowing-up with center $(x_s, x_{s+1})B$,
      hence a closed map.
    There exists a polynomial~$f$ over~$B$
      such that $\ideall_0 = \idealn C_0 + f C_0$
      and the leading coefficient of~$f$ is not contained in~$\idealn$.

    By Lemma~\ref{lem:2.2} and Theorem~\ref{thm:4.8},
      we have, for any $1 \leq k \leq s$,
        \begin{equation} \label{eqn:5.3.1}
          H_{\idealn B[T] + fB[T]}^p
          H_{(x_k, \dots, x_{s+1})}^q(B[T]) = 0
          \quad
          \text{if $p<d-s$ or $q=0$.}
        \end{equation}
    In fact,
      the leading coefficient of~$f$ is a regular element
      on $H_\idealn^{d-s} H_{(x_k, \dots, x_{s+1})}^q(B[T])$
      because it acts on
      the injective envelope of~$B/ \idealn$
      as isomorphism.
    Taking the local cohomology of
      a short exact sequence
        $$
          0 @>>>
          B[T] @>x_sT - x_{s+1}>>
          B[T] @>>>
          B[T]/ (x_sT - x_{s+1}) @>>> 0
        $$
      with respect to~%
        $
          (x_k, \dots, x_{s+1})
          = (x_k, \dots, x_s, x_sT - x_{s+1})
        $,
      we obtain an exact sequence
        \begin{multline*}
          0 @>>>
          H_{(x_k, \dots, x_{s+1})}^{s-k+1}(B[T]) @>>>
          H_{(x_k, \dots, x_{s+1})}^{s-k+1}(B[T]/(x_sT - x_{s+1})) @>>>
        \\
          @>>>
          H_{(x_k, \dots, x_{s+1})}^{s-k+2}(B[T]) @>>>
          H_{(x_k, \dots, x_{s+1})}^{s-k+2}(B[T]) @>>> 0,
        \end{multline*}
      because
      $(x_s, x_{s+1}) H_{(x_k, \dots, x_{s+1})}^{s-k+1}(B)=0$
      by Theorem~\ref{thm:4.8}.
    This and \eqref{eqn:5.3.1} show that
        $$
          H_{\idealn B[T] + fB[T]}^p H_{(x_k, \dots, x_{s+1})}^{s-k+1}
          (B[T]/ (x_s T - x_{s+1})) = 0
          \quad
          \text{for $p<d-s$}.
        $$
    Taking the local cohomology of an exact sequence
        $$
          0 @>>>
          \frac{(x_sT - x_{s+1}) \qtn \angled{x_s}}{(x_s T - x_{s+1})}
            @>>>
          B[T]/ (x_sT - x_{s+1}) @>>>
          C_0 @>>> 0
        $$
      with respect to $(x_k, \dots, x_{s+1})$,
      we obtain
        $$
          H_{(x_k, \dots, x_{s+1})}^{s-k+1}(C_0) =
          H_{(x_k, \dots, x_{s+1})}^{s-k+1}(B[T]/(x_s T - x_{s+1})),
        $$
      that is,
        \begin{equation} \label{eqn:5.3.2}
          H_{\ideall_0}^p H_{(x_k, \dots, x_{s+1})}^{s-k+1}(C_0) = 0
          \quad
          \text{for $p<d-s$}.
        \end{equation}

    We note that $x_s$ is $C_0$-regular.
    Put $k=s$.
    Then we have
        $$
          H_{\ideall_0}^p H_{(x_s, x_{s+1})}^q(C_0) = 0
          \quad
          \text{if $p<d-s$ or $q<1$}.
        $$
    By the spectral sequence
        $
          E_2^{pq} =
          H_{\ideall_0}^p H_{(x_s, x_{s+1})}^q(-) \Rightarrow
          H_{\ideall_0}^{p+q}(-),
        $
      we obtain
        \begin{equation} \label{eqn:5.3.3}
          H_{\ideall_0}^p(C_0) = 0
          \quad
          \text{for $p<d-s+1$,}
        \end{equation}
      that is,
      $\depth (C_0)_{\ideall_0} \geq d-s+1$.

    In the same way,
      we can show that
      $\depth (C_1)_{\ideall_1} \geq d-s+1$
      for any maximal ideal~$\ideall_1$
      of~$C_1$.
    Thus the first assertion is proved.
    In particular,
      $Y$ is Cohen-Macaulay if $s=1$.

    Assume that $A$ is Cohen-Macaulay.
    Using \cite[Lemma 1]{Faltings:78:Macaulay} twice,
      we find that
        $$
          x_{s+1} T_{s+2} - x_{s+2},
          \dots,
          x_{s+1} T_d - x_d,
          x_s T_{s+1} - x_{s+1}
        $$
      is a regular sequence on~$A[T_{s+1}, \dots, T_d]$.
    Therefore
        $$
          C_0 \cong
          A[T_{s+1}, \dots, T_d] /
          (x_{s+1} T_{s+2} - x_{s+2}, \dots, x_{s+1} T_d - x_d,
          x_s T_{s+1} - x_{s+1})
        $$
      is Cohen-Macaulay.
    In the same way, we can show that $C_1$ is Cohen-Macaulay.
    The proof is completed.
  \end{pf}

  In the rest of this section,
    we assume that $s \geq 2$ and
    let $\trans \idealb = \idealb \qtn \angled{x_{s-1}}$.

  \begin{lem} \label{lem:5.4}
    For any positive integer~$n$,
      $$
        \trans \idealb^n =
        \idealb^n \qtn \angled{x_{s-1}} =
        \idealq \idealb^{n-1} [(x_s, \dots, x_d) \qtn x_{s-1}]
        + x_s^n \idealq^{n-1} [\idealq \qtn x_{s-1}].
      $$
    In particular, $\trans \idealb^2 = \idealb \trans \idealb$.
  \end{lem}

  \begin{pf}
    It is sufficient to prove
        $$
          \idealb^n \qtn \angled{x_{s-1}} \subseteq
        \idealq \idealb^{n-1} [(x_s, \dots, x_d) \qtn x_{s-1}]
        + x_s^n \idealq^{n-1} [\idealq \qtn x_{s-1}].
        $$
    Take $a \in \idealb^n \qtn \angled{x_{s-1}}$.
    Then, by Lemma~\ref{lem:2.4}, Lemma~\ref{lem:3.2}
      and Theorem~\ref{thm:4.2},
      we have
        \begin{align*}
          a & \in (x_s, \dots, x_d)^{2n} \qtn \angled{x_{s-1}}
        \\
          & = (x_s, \dots, x_d)^{2n-1} [(x_s, \dots, x_d) \qtn x_{s-1}]
        \\
          & =
            [\idealq^{2n-1} + x_s \idealq^{2n-2} + \dots + (x_s^{2n-1})]
            [(x_s, \dots, x_d) \qtn x_{s-1}]
        \\
          & \subseteq \idealq \idealb^{n-1} [(x_s, \dots, x_d) \qtn x_{s-1}]
            + (x_s^n).
        \end{align*}
    If we put $a = b + x_s^n a'$
      where
        $b \in \idealq \idealb^{n-1} [(x_s, \dots ,x_d) \qtn x_{s-1}]$,
      then $x_s^n a' \in \idealb^n \qtn \angled{x_{s-1}}$.
    Since
      $x_{s-1}^l x_s^n a' \in \idealb^n$
      for a sufficiently large~$l$,
      we can put $x_{s-1}^l x_s^n a' = c + x_s^n d$
      where $c \in \idealq^{2n} + \dots + x_s^{n-1} \idealq^{n+1}$
      and $d \in \idealq^n$.
    Then
        $
          x_{s-1}^l a' - d \in
          \idealq^{n+1} \qtn \angled{x_s} =
          \idealq^n [\idealq \qtn x_s]
        $.
    Hence,
      $x_{s-1}^l a' \in \idealq^n$ and
        $
          a' \in
          \idealq^n \qtn \angled{x_{s-1}} =
          \idealq^{n-1} [\idealq \qtn x_{s-1}]
        $.
    The proof is completed.
  \end{pf}

  Therefore the Rees algebra~$\trans S = A[\trans \idealb t]$
    is finitely generated over~$S$.
  Let $\trans Y = \proj \trans S$.

  \begin{prop} \label{prop:5.5}
    $D_{(x_{s-1}, x_s, x_{s+1})}^0(S_+) = \trans S_+$.
  \end{prop}

  \begin{pf}
    First show that $x_{s-1}$, $x_s$ is
      an $\trans S_+$-regular sequence.
    Let $n > 0$.
    It is clear that $x_{s-1}$ is $\trans \idealb^n$-regular
      because it is $A$-regular.
    Let $a \in (x_{s-1} \trans \idealb^n \qtn x_s) \cap \trans \idealb^n$.
    Then $x_{s-1}^l a \in \idealb^n$
      for a sufficiently large~$l$.
    Since $x_s a \in (x_{s-1})$ and
      $x_s$,~\dots, $x_d$ is a d-sequence on~$A/ x_{s-1}^{l+1}A$,
        \begin{align*}
          x_{s-1}^l a & \in (x_{s-1}^{l+1}) \qtn x_s \cap \idealb^n
        \\
          & \subseteq (x_{s-1}^{l+1}) \qtn x_s \cap
            (x_{s-1}^{l+1}, x_s, \dots, x_d)
        \\
          & = (x_{s-1}^{l+1}).
        \end{align*}
    Hence $a \in (x_{s-1})$.
    If we put $a = x_{s-1} a'$,
      then
        $
          a' \in \idealb^n \qtn x_{s-1}^{l+1} \subseteq \trans \idealb^n
        $,
      that is, $a \in x_{s+1}\trans \idealb^n$.
    Thus we have proved that $x_s$ is
      $\trans S_+ / x_{s-1} \trans S_+$-regular.

    By~\eqref{eqn:2.1.1},
      we have
        \begin{equation} \label{eqn:5.5.1}
          D_{(x_{s-1}, x_s, x_{s+1})}^0(S_+) \subseteq
          D_{(x_{s-1}, x_s, x_{s+1})}^0(\trans S_+) = \trans S_+.
        \end{equation}
    Since $\idealq \qtn x_{s-1} \subseteq \idealq \qtn x_s$
      by Theorem~\ref{thm:4.2},
        $
          (x_{s-1}, \dots, x_d) \trans \idealb^n
          \subseteq \idealb^n
        $
      for all~$n>0$ by Lemma~\ref{lem:5.4},
      that is,
      $(x_{s-1}, \dots, x_d) \trans S_+ \subseteq S_+$.
    We have shown the inverse inclusion of~\eqref{eqn:5.5.1}.
  \end{pf}

  The following theorem is one of main aims of this section.

  \begin{thm}
    With notation as above,
        $$
          \depth \sheafO_{\trans Y, \trans q} \geq d-s+2
          \quad
          \text{for any closed point~$\trans q$ of~$\trans Y$}.
        $$
    In particular,
      if $s=2$,
      then $\trans Y$ is Cohen-Macaulay.
  \end{thm}

  \begin{pf}
    We have only to compute the depth of
        $$
          \trans C_0 = A[\trans \idealb/x_s x_{s+1}]
            \quad
            \text{and}
            \quad
          \trans C_1 = A[\trans \idealb/ x_{s+1}^2].
        $$
    Proposition~\ref{prop:5.5} says that
      $\trans C_i = D_{(x_{s-1}, x_s, x_{s+1})}^0(C_i)$
      and it is a finitely generated $C_i$-module
      for $i=0$,~$1$.

    Let $\trans\ideall_i$ be a maximal ideal of~$\trans C_i$
      and $\ideall_i = \trans\ideall_i \cap C_i$.
    Then $\ideall_i$ is a maximal ideal of~$C_i$
      because $\trans C_i$ is integral over~$C_i$.
    We use \eqref{eqn:5.3.2} as $k=s-1$,
      that is,
        \begin{equation} \label{eqn:5.6.1}
          H_{\ideall_i}^p H_{(x_{s-1}, x_s, x_{s+1})}^2(C_i) =0
          \quad
          \text{for $p<d-s$}.
        \end{equation}
    By using~\eqref{eqn:2.1.2},
      we obtain
        $$
          H_{\ideall_i}^p
            H_{(x_{s-1}, x_s, x_{s+1})}^q(\trans C_i) = 0
          \quad
          \text{if $p<d-s$ or $q<2$}.
        $$
    By the spectral sequence
        $
          E_2^{pq} =
          H_{\ideall_i}^p H_{(x_{s-1}, x_s, x_{s+1})}^q(-)
          \Rightarrow
          H_{\ideall_i}^{p+q}(-)
        $,
      we find
        \begin{equation} \label{eqn:5.6.2}
          H_{\ideall_i}^p(\trans C_i) = 0
          \quad
          \text{for $p<d-s+2$,}
        \end{equation}
      that is,
      $\depth (\trans C_i)_{\trans\ideall_i} \geq d-s+2$.
    Thus the proof is completed.
  \end{pf}

  The following corollary shall be used in the next section.

  \begin{cor} \label{cor:5.7}
    If $A/(x_s, \dots, x_d)$ is Cohen-Macaulay, then
        $$
          \depth \sheafO_{Y,q} \geq d-s+2
          \quad
          \text{for any closed point~$q$ of~$Y$.}
        $$
  \end{cor}

  \begin{pf}
    It is sufficient to prove $\trans \idealb=\idealb$.
    Let $a \in \trans \idealb$ and
      $l$ be an integer such that $x_{s-1}^l a \in \idealb$.
    Then we have
        \begin{align*}
          a & \in (x_s, \dots, x_d)^2 \qtn x_{s-1}^l
        \\
          & = (x_s, \dots, x_d) [(x_s, \dots, x_d) \qtn x_{s-1}^l]
        \\
          & = (x_s, \dots, x_d)^2 = \idealb + (x_s^2)
        \end{align*}
      by Lemma~\ref{lem:3.2}.
    Hence,
      we may assume that $a \in (x_s^2)$.
    Let $a = x_s^2 a'$.
    Since $x_{s-1}^l a \in \idealb \subseteq \idealq$,
      $a' \in \idealq \qtn x_{s-1}^l x_s^2 = \idealq \qtn x_s$
      by Theorem~\ref{thm:4.2}.
    Hence $a = x_s^2 a' \in x_s \idealq \subset \idealb$.
  \end{pf}

  We shall give another Macaulayfication of~$\spec A$
    by considering an ideal
      $\idealc = (x_{s-1}, \dots, x_d) \idealb$.
  Let $Z = \proj A[\idealc t]$,
    which is the blowing-up of~$Y$
    with center $(x_{s-1}, \dots, x_d) \sheafO_Y$.

  \begin{thm} \label{thm:5.8}
    With notation as above,
        $$
          \depth \sheafO_{Z,r} \geq d-s+2
          \quad
          \text{for any closed point~$r$ of~$Z$}.
        $$
    Furthermore,
      if $s=2$ or $A$ is Cohen-Macaulay,
      then $Z$ is Cohen-Macaulay.
  \end{thm}

  \begin{pf}
    Since
        $
          (x_{s-1}x_s, x_s^2) \idealq +
          x_{s-1} (x_{s+1}^2, \dots, x_d^2) +
          (x_{s+1}^3, \dots, x_d^3)
        $
      is a reduction of~$\idealc$,
      we have only to compute the depth of
        \begin{align*}
          D_0 & = A[\idealc/ x_{s-1} x_s x_{s+1}] = C_0[x_s/ x_{s-1}],
        \\
          D_1 & = A[\idealc/ x_s^2 x_{s+1}] = C_0[x_{s-1}/ x_s],
        \\
          D_2 & = A[\idealc/ x_{s-1} x_{s+1}^2] = C_1[x_{s+1}/ x_{s-1}],
        \\
          \intertext{and}
          D_3 & = A[\idealc/ x_{s+1}^3] = C_1[x_{s-1}/ x_{s+1}].
        \end{align*}
    For $i=0$ or~$1$,
      let $\ideall_i$ be a maximal ideal of~$C_i$.
    By~\eqref{eqn:2.1.1},
      there exists an exact sequence
        $$
          0 @>>>
          C_i @>>>
          \trans C_i @>>>
          H_{(x_{s-1}, x_s, x_{s+1})}^1(C_i) @>>> 0.
        $$
    By using~\eqref{eqn:5.3.3}
      and \eqref{eqn:5.6.2},
      we obtain
        $$
          H_{\ideall_i}^p
          H_{(x_{s-1}, x_s, x_{s+1})}^1(C_i) = 0
          \quad
          \text{for all $p<d-s$.}
        $$
    Furthermore,
      $(x_{s-1}, \dots, x_d) \trans C_i \subseteq C_i$:
      see the proof of Proposition~\ref{prop:5.5}.
    Therefore, by~\eqref{eqn:5.6.1}, we have
        \begin{equation} \label{eqn:5.8.1}
          H_{\ideall_i}^p H_{(x_{s-1}, x_s, x_{s+1})}^q(C_i) = 0
          \quad
          \text{if $p<d-s$ or $q = 0$}
        \end{equation}
      and
        \begin{equation}
            \label{eqn:5.8.2}
          (x_{s-1}, \dots, x_d)
          H_{(x_{s-1}, x_s, x_{s+1})}^1(C_i) =0.
        \end{equation}
    Therefore we can prove
        $$
          \depth (D_i)_{\idealr_i} \geq d-s+2
        $$
      for any maximal ideal~$\idealr_i$
      of~$D_i$ and $i=0$,~\dots, $3$
      in the same way as Theorem~\ref{thm:5.3}.

    To make sure,
      we compute the depth of~%
        $
          D_0 \cong C_0[T]/(x_{s-1}T- x_s) \qtn \angled{x_{s-1}}
        $.
    First we note that $x_{s+1} \in x_s C_0$
      and $x_{s+1} \in x_s D_0$.
    Let $\idealr_0$ be a maximal ideal of~$D_0$
      and $\ideall_0 = \idealr_0 \cap C_0$.
    Then $\ideall_0$ is a maximal ideal of~$C_0$
      and there exists a polynomial~$f$
      over~$C_0$
      such that $\idealr_0 = \ideall_0 D_0 + f D_0$
      and the leading coefficient of~$f$
      is not contained in~$\ideall_0$.
    We obtain
        $$
          H_{\ideall_0 C_0[T] + f C_0[T]}^p
          H_{(x_{s-1}, x_s)}^q (C_0[T]) = 0
          \quad
          \text{if $p<d-s+1$ or $q = 0$}
        $$
      from~\eqref{eqn:5.8.1}.
    Taking the local cohomology of an exact sequence
        $$
          0 @>>>
          C_0[T] @>x_{s-1}T - x_s>>
          C_0[T] @>>>
          C_0[T]/(x_{s-1}T - x_s) @>>>
          0,
        $$
      we have an exact sequence
        \begin{multline*}
          0 @>>>
          H_{(x_{s-1}, x_s)}^1(C_0[T]) @>>>
          H_{(x_{s-1}, x_s)}^1(C_0[T]/ (x_{s-1}T - x_s)) @>>>
        \\
          @>>>
          H_{(x_{s-1}, x_s)}^2(C_0[T]) @>>>
          H_{(x_{s-1}, x_s)}^2(C_0[T]) @>>>
          0
        \end{multline*}
      because of~\eqref{eqn:5.8.2}.
    This says that
        $$
          H_{\ideall_0 C_0[T] + f C_0[T]}^p
          H_{(x_{s-1}, x_s)}^1
            (C_0[T]/ (x_{s-1}T - x_s)) = 0
          \quad
          \text{for $p<d-s+1$}.
        $$
    Taking the local cohomology of an exact sequence
        $$
          0 @>>>
          \frac{(x_{s-1}T-x_s) \qtn \angled{x_{s-1}}}
            {(x_{s-1}T - x_s)} @>>>
          C_0[T]/(x_{s-1}T - x_s) @>>>
          D_0 @>>>
          0
        $$
      with respect to~$(x_{s-1}, x_s)$,
      we obtain
        $$
          H_{\idealr_0}^p H_{(x_{s-1}, x_s)}^1(D_0) = 0
          \quad
          \text{for $p<d-s+1$}.
        $$
    Of course,
      $H_{(x_{s-1}, x_s)}^0(D_0) = 0$.
    By the spectral sequence
        $$
          E_2^{pq} =
          H_{\idealr_0}^p H_{(x_{s-1}, x_s)}^q(-) \Rightarrow
          H_{\idealr_0}^{p+q}(-),
        $$
      we get $H_{\idealr_0}^p(D_0) = 0$
      for any $p<d-s+2$.
    That is, $\depth (D_0)_{\idealr_0} \geq d-s+2$.

    The last assertion is also proved
      in the same way as Theorem~\ref{thm:5.3}.
    \end{pf}

\setcounter{equation}{0}
\section{The proof of Theorem~\ref{mthm}}
  \label{sec:6}

  This section is devoted to the proof of Theorem~\ref{mthm}.
  Let $A$ be a Noetherian ring possessing a dualizing complex
    and $X$ a quasi-projective scheme over~$A$.
  That is, $X$ is a dense open subscheme
    of $X\closure = \proj R$
    where $R = \bigoplus_{n \geq 0} R_n$ is a Noetherian graded ring
    such that $R_0$ is a homomorphic image of~$A$ and
    $R$ is generated by~$R_1$ as an $R_0$-algebra.
  Let $V\closure$ be the non-Cohen-Macaulay locus of~$X\closure$
    and $U\closure = X\closure \setminus V\closure$.
  Of course $V = V\closure \cap X$ is
    the non-Cohen-Macaulay locus of~$X$.
  Let $\dc$ be a dualizing complex of~$R$
    with codimension function~$v$.
  Assume that $X$ satisfies the assumption of Theorem~\ref{mthm}.

  Without loss of generality,
    we may assume that
      \begin{equation} \label{eqn:6.0.1}
        v(\idealp) = 0
          \quad
        \text{for all associated prime ideal~$\idealp$ of~$R$:}
      \end{equation}
    see \cite[p.~191]{Faltings:78:Macaulay}.
  Then the local ring~$\sheafO_{X,p}$ of~$p \in X$ satisfies
    the assumption of Section~\ref{sec:5},
    that is,
    $\dim \sheafO_{X,p}/ \idealp = \dim \sheafO_{X,p}$
    for any associated prime ideal~$\idealp$ of~$\sheafO_{X,p}$.
  For the sake of completeness,
    we sketch out the proof.
  Let $\ideala$ be a homogeneous ideal of~$R$
    such that $V\closure = V(\ideala)$.
  Then the closed immersion $\proj R/H_\ideala^0(R) \to X\closure$
    is birational as follows.
  For any minimal prime ideal~$\idealp$ of~$R$,
    $\ideala \not\subset \idealp$ and
    $H_\ideala^0(R) \subseteq \idealp$
    because $R_\idealp$ is Cohen-Macaulay.
  Hence the underlying set of~$\proj R/H_\ideala^0(R)$
    coincides with the one of~$X\closure$.
  Furthermore,
    $f^{-1}(U\closure) \to U\closure$ is an isomorphism and
    $U\closure$ is dense in~$X\closure$.
  By replacing $R$ by~$R/H_\ideala^0(R)$,
    we may assume that
      \begin{equation} \label{eqn:6.0.2}
        \text{
          every associated prime ideal of~$R$
          is minimal.
        }
      \end{equation}
    Next we fix a primary decomposition of~$(0)$ in~$R$.
  For all integer~$i$,
    let $\idealq_i$ be the intersection of
    all primary component~$\idealq$ of~$(0)$
    such that $v(\sqrt\idealq) = i$.
  Then $g \colon \coprod_i \proj R/\idealq_i \to X\closure$
    is a finite morphism
    and $g^{-1}(U\closure) \to U\closure$ is an isomorphism
    as follows.
  Note that $\idealq_i = R$ for all but finitely many~$i$.
  Furthermore, for any $\idealp \in U\closure$,
    $\idealp \supseteq \idealq_i$ if and only if
    $v(\idealp) - \dim R_\idealp = i$
    because $R_\idealp$ is Cohen-Macaulay, hence equidimensional.
  Therefore $U\closure$ is
    the disjoint union of $U\closure \cap V(\idealq_i)$.
  Moreover $R\hlz\idealp = [R/\idealq_i]\hlz\idealp$
    if $\idealp \in U\closure \cap V(\idealq_i)$.
  Because of~\eqref{eqn:6.0.2},
    $g^{-1}(U\closure)$ and $U\closure$ are dense
    in $\proj R/\idealq_i$ and $X\closure$,
    respectively.
  Thus $g^{-1}(X) \to X$ is birational proper
    and the connected components of~$g^{-1}(X)$
    satisfy the assumption of Theorem~\ref{mthm}.

  Since $u$ is locally constant,
    $V_i = u^{-1}(i) \cap V$ is
    closed for any positive integer~$i$.
  We put $s_i = \dim V_i$.
  By~\eqref{eqn:6.0.1},
    we find that $V_1 = \emptyset$, $s_2 \leq 0$ and $s_3 \leq 1$.
  Let $d$ be the largest integer such that $V_d \ne \emptyset$
    and $s = s_d$.
  We shall give a closed subscheme~$W$ of~$X$
    such that $V_d = V \cap W$ and
    $\sheafO_{Y,q}$ is Cohen-Macaulay
    for all $q \in \pi^{-1}(W)$
    where $\pi \colon Y \to X$ is the blowing-up
    of~$X$ with center~$W$.
  Let $\ideala = \prod_{i>0} \ann H^i(\dc)$,
    which is finite product.
  Then it is obvious that $V\closure = V(\ideala)$.
  Fix a primary decomposition of~$\ideala$
    and let $\ideala_d$ be the intersection of
    all primary component~$\idealq$ of~$\ideala$
    such that $\sqrt\idealq \in V_d$.
  Then we can take homogeneous elements $z_1$,~\dots, $z_d \in R$
    such that
      \begin{gather} \label{eqn:6.0.3}
        V_i \cap V((z_{d-s_i}, \dots, z_d)) = \emptyset
        \quad \text{for $i < d$}
      \\
        d(\idealp) = d
        \quad
        \text{for all minimal prime ideal~$\idealp$
          of $R/(z_1, \dots, z_d) \qtn \angled{R_+}$}
      \\
        \begin{cases}
          z_{s+1}, \dots, z_d \in \ideala_d;
        \\
          \label{eqn:6.0.5}
          z_i \in
          \prod_{j>d-i} \ann H^j(\hom(R/(z_{i+1}, \dots, z_d), \dc)), &
          \text{for $i \leq s$}
        \end{cases}
      \end{gather}
    in the same way as Section~\ref{sec:4}.
  We put
      $$
        \idealb =
        \begin{cases}
          (z_1, \dots, z_d), & \text{if $s=0$};
        \\
          (z_1, \dots, z_d) (z_2, \dots, z_d), & \text{if $s=1$};
        \\
          (z_1, \dots, z_d) (z_2, \dots, z_d) (z_3, \dots, z_d),
          & \text{if $s=2$}
        \end{cases}
      $$
    and prove that $W = V(\idealb) \cap X$ satisfies
    the required properties.

  Because of~\eqref{eqn:6.0.3},
    $V_i \cap W = \emptyset$ for~$i < d$.
  Let $\pi \colon Y \to X$ be the blowing-up of~$X$
    with center~$W$,
    $q$ a closed point of~$\pi^{-1}(W)$ and
    $\idealp \subseteq R$ the image of~$q$.
  Take an element $y \in R_1 \setminus \idealp$
    and put $x_i = z_i / y^{\deg z_i}$ for all~$i$.
  Since $(\dc)\hlz\idealp$ is a dualizing complex of~$R\hlz\idealp$,
    we obtain
      $$
        \begin{cases}
          x_{s+1}, \dots, x_d \in \ideala_{R\hlz\idealp}(R\hlz\idealp);
        \\
          x_i \in \ideala_{R\hlz\idealp}(
            R\hlz\idealp/(x_{i+1}, \dots, x_d)),
          & \text{for $i \leq s$}.
        \end{cases}
      $$
    from~\eqref{eqn:6.0.5}.

  When $s=2$, there exist three cases:
  If $z_1$, $z_2 \in \idealp$,
    then $x_1$,~\dots, $x_d$ is a system of parameters
    for~$R\hlz\idealp$ satisfying~\eqref{eqn:4.0.1}
    or a regular sequence
    on the Cohen-Macaulay ring~$R\hlz\idealp$.
  Since
      $
        \idealb\hlz\idealp =
        (x_1, \dots, x_d) (x_2, \dots, x_d) (x_3, \dots, x_d)
      $,
    $\sheafO_{Y,q}$ is Cohen-Macaulay by Theorem~\ref{thm:5.8}.

  If $z_2 \in \idealp$ but $z_1 \notin \idealp$,
    then $x_2$,~\dots, $x_d$ is a subsystem of parameters
    for~$R\hlz\idealp$ satisfying~\eqref{eqn:4.0.1}
    or a regular sequence
    on the Cohen-Macaulay ring~$R\hlz\idealp$.
  Furthermore
    $\idealb\hlz\idealp = (x_2, \dots, x_d)(x_2, \dots, x_d)$
    and $R\hlz\idealp/(x_2, \dots, x_d)$ is Cohen-Macaulay
    because
    $x_1 \in \ideala_{R\hlz\idealp}(R\hlz\idealp/(x_2, \dots, x_d))$
    is a unit.
  Hence $\sheafO_{Y,q}$ is Cohen-Macaulay by Corollary~\ref{cor:5.7}.

  If $z_1$, $z_2 \notin \idealp$, then
    $x_3$,~\dots, $x_d \in \ideala_{R\hlz\idealp}(R\hlz\idealp)$
    is a subsystem of parameters for~$R\hlz\idealp$
    and $R\hlz\idealp/(x_3, \dots, x_d)$ is Cohen-Macaulay.
  Since $\idealb\hlz\idealp = (x_3, \dots, x_d)$,
    $\sheafO_{Y,q}$ is Cohen-Macaulay by Theorem~\ref{thm:5.1}.

  When $s=0$ or $1$,
    we can prove the assertion
    in the same way as above.

  By repeating this procedure,
    we obtain a Macaulayfication of~$X$.
  We complete the proof of Theorem~\ref{mthm}.

\subsection*{Acknowledgment}
  The author is grateful
    to Professor S.~Goto
    for his helpful discussions
    and to Professor K.~Kurano
    for careful reading of the draft.

\ifx\undefined\bysame
  \newcommand{\bysame}{%
    \leavevmode\hbox to3em{\hrulefill}\,}
\fi

\enddocument